**Bootstrap Markov chain Monte Carlo and optimal solutions for the Law of**

**Categorical Judgment (Corrected)**


Greg Kochanski and Burton S. Rosner

University of Oxford



Corresponding author: Burton S Rosner, 2 Ryeford Court, Wish Ward, Rye

TN31 7DH, U. K.; Tel.: (44)1797225714; email: **burton.rosner@phon.ox.ac.uk**.






Abstract

A novel procedure is described for accelerating the convergence of Markov chain Monte Carlo computations. The algorithm uses an adaptive bootstrap technique to generate candidate steps in the Markov Chain. It is efficient for symmetric, convex probability distributions, similar to multivariate Gaussians, and it can be used for Bayesian estimation or for obtaining maximum likelihood solutions with confidence limits.  As a test case, the Law of Categorical Judgment (Corrected) was fitted with the algorithm to data sets from simulated rating scale experiments.  The correct parameters were recovered from practical-sized data sets simulated for Full Signal Detection Theory and its special cases of standard Signal Detection Theory and Complementary Signal Detection Theory.

The Markov-Chain Monte-Carlo (MCMC) technique originated in statistical physics (Metropolis, Rosenbluth, Rosenbluth, Teller, & Teller, 1953). It has since spread into other fields, including psychology (e. g., Béguin & Glas, 2001; Griffiths, Steyvers,  & Tenenbaum, 2007; Morey, Rouder, & Speckman, 2008; Sanborn, Griffiths, & Shiffrin, 2010). Typically, MCMC is used to compute samples from the stationary probability distribution $\pi$ of a Markov chain. In the MCMC algorithm, only ratios of $\pi$ (or of likelihoods) need be calculated, so normalization of $\pi$ becomes unnecessary. Eliminating this normalization is a major advantage, since it requires a multidimensional integral over $\pi$ that is often computationally impractical.

The original MCMC algorithms (e.g., Metropolis et al., 1953; Hastings, 1970) were simple. From position $z_i$, compute a test position $x_i = z_i + v$, where $v$ is drawn from a distribution $V$ that generates the steps of a random walk. Then accept $x_i$ with probability $\min[\pi(x_i)/\pi(z_i); 1]$, either setting $z_{i+1} = x_i$ or else rejecting $x_i$ and setting $z_{i+1} = z_i$.

The generator $V$ is relatively unimportant asymptotically, because almost any reasonable $V$ will make the MCMC chain converge on the same distribution $\pi$. Practically, however, $V$ is critical, especially for high-dimensional problems. Bad choices of $V$ can dramatically increase the time needed to sample all of $\pi$. To address this problem, procedures have been devised that compute $V$ based on the history of the chain (e.g., Gelfand & Sahu, 1994; Gilks, Roberts, & Sahu, 1998; Atchadé & Rosenthal, 2005).



We now describe a novel algorithm, termed bootstrap Markov chain Monte Carlo (BMCMC), that also uses the chain's history to construct $V$. In the latter half of this paper, BMCMC is applied to obtain maximum likelihood solutions to a general equation for rating data, the Law of Categorical Judgment (Corrected) (Rosner & Kochanski, 2009). Our procedures go well beyond the earlier treatment of rating data by Schönemann and Tucker (1967).

## Bootstrap Markov Chain Monte Carlo

Bootstrap Markov chain Monte Carlo rests on standard Markov chain Monte Carlo procedures (Anonymous, 1998; Geyer,1992; Gilks, Richardson, & Spiegelhalter, 1995). Like any MCMC algorithm, it takes iterative steps, constantly updating a vector of parameters $\vec{p}$. The algorithm, however, adapts itself to the probability distribution. It adjusts its steps to achieve efficient sampling of probability densities that are approximately multivariate-Gaussian; but it proceeds in such a way that violations of the Markov assumptions become asymptotically negligible. It can produce chains of samples that converge rapidly to $\pi$ for approximately symmetric and convex densities.

The algorithm is robust. It can handle functions that are not computable (in the practical sense) at certain points; and it does not assume that $\pi$ is smooth or continuous, so long as it is bounded above. Finally, BMCMC can be used on functions that stochastically approximate $\pi$, where the computable approximation to $\pi$ is a random variable, and successive computations of $\pi$ with the same parameters may



give different values.[1] Experience has shown that this makes it workable for many types of computation where common procedures are unsuitable.  One example is computation of $\pi$ by sampling techniques such as a Monte-Carlo integral.  Another example is where $\pi$ is a noisy measurement of some physical property.

The BMCMC procedure can be run in either of two closely related modes: optimisation or sampling. In optimisation mode, it searches for parameter values that best explain the data under analysis, by maximizing log likelihood, perhaps modified by a Bayesian prior. (We confine attention here to the widely used log$L$, but this mode can also be used for Bayesian model fitting.) In sampling mode, BMCMC repeatedly samples from $\pi$ (proportional to log$L$ or to a Bayesian posterior density). This allows calculation of all necessary statistics, including confidence limits for the optimal parameters and for the observations predicted from those parameters.

The following description of BMCMC is somewhat simplified. The source code should be consulted for details; footnotes give references into the code. The algorithm is embodied in two main python modules, mcmc.py and mcmc_helper.py. They are in release gmisclib-0.65.5 which, along with all other code described here, can be downloaded from http://sourceforge.net/projects/speechresearch or under http://phon.ox.ac.uk.

Braun, Kochanski, Grabe, and Rosner (2006) and Alvey, Orphanidou, Coleman, McIntyre, Golding, and Kochanski (2008) used earlier versions of this code (albeit with minimal description). In other work, we have successfully tested BMCMC on moderately high-dimensional problems. We also have shown



convergence to expected solutions and reasonable error bars for 200-dimensional problems.

Application of BMCMC requires a user-supplied software module that produces the log of a value proportional to $\pi$. In optimisation mode, for example, $\pi(\vec{p})$ would be the probability with which a model with parameters $\vec{p}$ would generate the observed data. The problem-specific module takes a vector of tentative parameters $\vec{p}$ from the BMCMC process. The module first may modify $\vec{p}$ to account for symmetries or constraints on the distribution.[2] Then it computes the log likelihood $\log L$ for the resulting parameters,[3] finally returning $\log L$ to the BMCMC algorithm.

*Optimisation Mode*

The central operation of the BMCMC algorithm increments the current parameter vector $\vec{p}$ by a quasi-random step vector $\vec{d}$. The difference $D = \log L(\vec{p} + \vec{d}) - \log L(\vec{p})$ is computed. If $\log L$ improves ($D > 0$), the new location is accepted. If $D$ is substantially negative, or $\log L(\vec{p} + \vec{d})$ cannot be computed, the step is rejected, leaving the procedure at its prior location. For small decreases, the new location is randomly accepted or rejected. Then the process iterates from its current location. As the algorithm progresses, $\log L$ generally increases, and an optimal solution for the parameters finally emerges.

*Step determination*. Bootstrap Markov chain Monte Carlo uses two schemes for generating $\vec{d}$. After an initial start-up period, an adaptive bootstrap resampling



procedure is used to generate 90 per cent of the steps.[4] Otherwise, the step comes

from a pre-specified multivariate Gaussian density (with adaptive, diagonal scaling).[5]

The bootstrap procedure randomly chooses two vectors from an archive of

previously accepted parameter vectors. The difference between them becomes the

current step. In principle, bootstrapping violates the Markov assumption that each

step is independent of its predecessors, because the archive contains the history of the

algorithm's computations. However, as the archive lengthens, the density of samples

from it asymptotically approaches $\pi$ and becomes stationary.  Therefore, when the

archive holds a sufficiently large number of samples widely distributed across $\pi$, the

Markov assumption will be satisfied to any required accuracy.  We later show that the

archive size becomes asymptotically infinite, making the system approach a Markov

chain.

The bootstrapping scheme works well when the density of $\log L$ is close to a

multivariate Gaussian, even a highly elongated one. This is because a large enough

archive makes the step probability density approach the convolution

$$P(\vec{d}) = \int \pi(\vec{p}) * \pi(\vec{p} + \vec{d}) d\vec{p}.$$ If $\pi(\vec{p})$ is a multivariate Gaussian, then $P(\vec{d})$ will also

be. It will have the same shape and orientation as $\pi(\vec{p})$ but twice the variance. As a

result, the long axes of the distribution will be accurately aligned with those of $\pi$.

The BMCMC algorithm tracks the proportion $f_\lambda$ of bootstrap steps accepted

over a period $t_\lambda$. This proportion controls a factor $\lambda$ that multiplies $\vec{d}$ . When $f_\lambda$

significantly exceeds 1/4, $\lambda$ is increased, thus decreasing the acceptance rate. When



the proportion falls significantly below 1/4, $\lambda$ is decreased with the opposite effect.[6] If $\pi(\vec{p})$ is multivariate Gaussian, $\lambda$ converges to its final value within a few hundred steps. Its final value depends on the number of dimensions and on the details of $\pi(\vec{p})$ but is typically slightly smaller than unity.

This procedure for controlling $\lambda$ can sometimes behave badly if $\pi(\vec{p})$ drops rapidly, especially discontinuously, to zero, e.g. in an optimization problem where $\vec{p}$ is subject to hard-wall constraints. If three or more hard walls meet, they can form a corner where locally even small steps are less than 25 per cent likely to be accepted. This lowers the average step acceptance rate, and the algorithm responds by reducing $\lambda$, which in turn slows computation. If hard-wall constraints are desired, see the Appendix.

In principle, adjusting $f_\lambda$ also violates the Markov assumption, because the scaling depends on recent history. The period $t_\lambda$, however, increases during optimisation. Eventually, $t_\lambda$ becomes longer than the time required for BMCMC to explore all of $\pi$. The dependence of $\lambda$ on history then becomes unimportant, and the algorithm asymptotically produces a Markov chain. We insure that this happens by making the required significance level for deciding when to change $f_\lambda$ an increasingly stringent function of the number of iterations since the most recent reset.[7] Resets are described below.

The pre-specified step generator draws from a multivariate Gaussian density $V*$ for the first few steps. Then a scaled version of $V*$ is used to generate 10 per cent



of the later steps. The scale factor similarly depends on the fraction of recently accepted steps and is proportional to the square root of the parameter-by-parameter standard deviation of the archive.[8] Again, as $t_\lambda$ and the archive lengthen, this step generator also asymptotically behaves as a Markov chain.

The two methods of step selection compensate for each other's deficiencies. Alone, the pre-specified step distribution can cause slow convergence and slow exploration of $\pi(\vec{p})$ if the shape or orientation of $V^*$ does not match that of $\pi$. Nevertheless, this method will eventually explore all of $\pi(\vec{p})$.

The bootstrap method needs the pre-specified distribution procedure to initialise the archive. Furthermore, bootstrapping is limited to linear combinations of archive points. Hence, if all archived vectors fell in a low dimensional subspace $\varpi$ of $\pi(\vec{p})$, bootstrapping would remain there. Mixing the two methods of step selection avoids that trap, because picking from the pre-specified distribution will soon engineer an escape from $\varpi$.

*Operation*. Optimisation mode is broadly similar to simulated annealing (Press, Teukolsky, Vettering, & Flannery, 2002, pp. 448-458), including a system temperature.[9] Early on, large decreases in $\pi$ are allowed, corresponding to a high system temperature. The default annealing schedule decreases the temperature whenever a step is accepted, eventually approaching a specified target temperature.[10] (The annealing schedule can be redesigned by changing parameters or re-implementing the mcmc_helper.step_acceptor object.) The algorithm stores the



maximum value of log$L$ encountered so far.  If a new value appears that exceeds the old one by half the temperature, the algorithm is partially reset: the temperature is raised slightly, all counters are reset, $t_\lambda$ is decreased, and the archive is shortened by eliminating those $\vec{p}$ that give the smallest log$L$.[11] Additionally, over the next 2*$N_p$ accepted steps, each step begins at the parameters that give the maximal log$L$ rather than the most recent parameter vector.[12]  (Here, $N_p$ is the number of unconstrained parameters.) This can substantially improve the rate at which the algorithm first approaches a minimum. The search then continues.

Resets serve several purposes. Raising the temperature can help the algorithm to escape a local maximum. Shortening the archive and $t_\lambda$ facilitates adaptation of the search to the shape of log$L$.

For any $\pi$ with an upper bound, the tolerance of 0.5$T$ for log$L$ guarantees that there will be only a finite number of resets. Consequently, there will be a final reset, after which both the archive size and the tracking period $t_\lambda$ will approach infinity. This is necessary for BMCMC to become asymptotically Markovian.

*Optimisation termination.*  The BMCMC optimisation mode terminates when two conditions are met.[13] First, successive accepted values of $\vec{p}$ must show no systematic drift. Second, enough steps must have occurred to match an estimate of the number needed to explore all of $\pi$.

Drift is evaluated by sampling archived $\vec{p}$ every $N_a$ accepted steps ($N_a$ is set at 24 in the examples below). Difference vectors are computed between adjacent



samples, and the angles between pairs of successive difference vectors are found. If $\vec{p}$ systematically drifts, the difference vectors will point in a common direction. Otherwise, reversals in direction will occur, and the angles will frequently exceed $\pi/2$ radians. The algorithm counts the number of such reversals since the last reset.[14] The termination condition of no systematic drift is met when the count reaches a suitable threshold.

The second termination condition requires an estimate of the number of steps needed to explore all of $\pi$. The algorithm uses a bent multivariate Gaussian model to make this estimate. The density of $\log L$ is represented as approximately multivariate Gaussian, but the longest axis of the probability ellipsoid is assumed to be slightly curved. This curvature limits the length of steps along the bent direction. Since steps are straight, long steps would fall off the likelihood ridge, yielding small values of $\log L$ that typically would be rejected. These small steps then form a random walk along the longest axis.

The curvature in the model is related to the factor $\lambda$ that scales the bootstrap step size. It should take on the order of $\lambda^{-2}$ accepted steps for a random walk to explore the length of the density.[15] Accordingly, the second termination condition is that the number of steps accepted since the last reset exceeds a constant times $\lambda^{-2}$.[16] (N.B.: if there are $k \geq 2$ comparably large eigenvalues of the covariance matrix, this termination condition may fire early.)



*Sampling Mode*

After optimisation finishes, the BMCMC procedure can be run in sampling mode, producing various confidence intervals.[17] The algorithm randomly samples a probability density $\pi(\vec{p})$, proportional to $\log L(\vec{p})$, in the vicinity of the solution. The system temperature is held at unity, and steps are accepted in accordance with the Metropolis algorithm. Each sample of $\vec{p}$ is archived.

Sampling mode sometimes finds a new maximum of $\log L$. If so, the previously described reset procedures apply, except that now the oldest archived samples are dropped and the temperature remains fixed at unity. In sampling mode, the archive becomes very long after the final reset. If it needs shortening, e.g., to limit computer memory consumption, the oldest samples would be dropped. Given sufficient time and memory, sampling from the archive becomes a good approximation of sampling from $\pi(\vec{p})$.

*Sampling termination.* The BMCMC algorithm stops when the desired number of samples, $N_e$, from $\pi(\vec{p})$ have been stored and enough steps have been accepted to explore the longest axis of the density approximately $N_e$ times.[18] Iterations continue as long as $\sum_i 1 / \lambda_i^2 < 1.4 N_p N_e$, where $\lambda_i$ is the step scaling factor $\lambda$ at the $i$th iteration and $N_p$ is the number of free parameters. This proviso flows from the same bent random walk model that imposes a termination condition on optimisation mode. Within the limits of that model, the inequality insures getting least $N_e$ independent samples from $\pi(\vec{p})$. The summation starts anew whenever the



algorithm resets itself, so iterations cannot terminate while log$L$ continues to improve.

     *Confidence intervals*. The archive resulting from repeated sampling yields confidence limits for the parameters and for the predicted data. All the usual statistics, such as the means and variance of parameters, are integrals over $\pi(\vec{p})$. Any such integral can be approximated by a sum over samples from $\pi(\vec{p})$. Confidence intervals can be computed from estimates of means and variances or determined directly from cumulants of the samples.

<div align="center">Optimal Solutions for the Law of Categorical Judgment (Corrected)</div>

     We now illustrate the use of BMCMC for multidimensional nonlinear optimisations, including computation of error bars. The Law of Categorical Judgment (Torgerson, 1958; McNicol, 1972) is a general model for any rating experiment. On each trial $t$, one stimulus $S_h$ is randomly selected for presentation from a pool of $N$ stimuli. The subject responds with one of $M+1$ possible ordered responses $R_i$.

     According to the law, each stimulus projects an independent random Gaussian density $\varphi(x;\mu_{S_h},\sigma_{S_h})$ of possible impressions on a one-dimensional psychological continuum $v$. On trial $t$, the observer draws one sample with value $s_t$ on $v$ from the signal density for $S_h$. Only the experimenter knows which $S_h$ actually produced $s_t$. On each trial, the observer also places $M$ criteria at different loci on $v$. The criteria divide $v$ into $M+1$ successive intervals that define the possible responses.



The Law of Categorical Judgment proposes that each criterion also is an independent random Gaussian variable, $\varphi(x; \mu_{C_i}, \sigma_{C_i})$, on $v$. On trial $t$, therefore, criterion positions $c_{it}$, $i = 1, \ldots, M$, are drawn from the criterion densities, where $i$ indexes the order of the density means. The indices are private to the observer. The observer calculates the $M$ differences $s_t$-$c_{it}$. If $s_t$-$c_{jt}$ is the smallest positive difference, response $R_j$ is selected. The response $R_{M+1}$ ensues when all differences are negative.

*The Law of Categorical Judgment (Corrected)*

The long-accepted general equation for the Law of Categorical Judgment (Torgerson, 1958; McNicol, 1972) is fatally flawed (Rosner & Kochanski, 2009). It ignores the fact that criterion densities are independent Gaussians. When the densities overlap sufficiently, criterion samples can get out of order, and the original formula can yield negative probabilities. Rosner and Kochanski derived a new equation, the Law of Categorical Judgment (Corrected), that fixes the flaw:

$$P(\text{R} = i | \text{S}_h) = \int_{-\infty}^{\infty} \varphi(c_i; \mu_{C_i}, \sigma_{C_i}) \int_{-\infty}^{c_i} \varphi(s; \mu_{S_h}, \sigma_{S_h}) \prod_{j \neq i}^{M} [1 - \int_{s}^{c_i} \varphi(c_j; \mu_{C_j}, \sigma_{C_j}) dc_j] ds dc_i. \qquad (1)$$

Equation 1 assumes that criteria are independent Gaussian variables. Nevertheless, it still may apply when response dependencies undercut the assumption. Criterion-setting theory (Treisman & Williams, 1984) explains how such dependencies arise in detection and discrimination. Treisman (1985) extended the theory to absolute judgments, a form of rating. Simulations showed that the resulting



criterion densities were indistinguishable from Gaussians. In effect, response dependencies would drive the product term in Equation 1 towards unity. The product term, however, would not necessarily reach unity. When Treisman and Faulkner (1985) analysed rating data with criterion-setting theory, they found evidence that criteria interchanged positions. Equation 1 therefore could apply to rating data when response dependencies occurred.

*Two special cases.* The familiar signal detection theory (SDT) model for rating experiments makes criterion positions $c_{it}$ invariant across trials (see Macmillan & Creelman, 2003, ch. 3). Each criterion density thereby becomes a Dirac delta function $\delta(x-c_i)$ (Bracewell, 1999, chap. 5), and the outer integral and product in Equation 1 simplify away. The result (Rosner and Kochanski, 2009) is

$$P(\text{R}=i|\text{S}) = \int_{c_{i-1}}^{c_i} \varphi(s;\mu_\text{S},\sigma_\text{S})ds. \tag{2}$$

This is equivalent to the standard form for the SDT rating model:

$$z(\text{R} \leq i \mid \text{S}_h) = (\mu_{C_i} - \mu_{S_h})/\sigma_{S_h}. \tag{3}$$

To emphasize its relationship to the SDT model (Equations 2 and 3), Equation 1 is referred to as Full Signal Detection Theory (FSDT).

In the limit opposite to SDT, all criterion densities are Gaussian but all signal densities Dirac delta functions $\delta(x-s)$. The equation for this complementary signal



detection (CSDT) model (Rosner & Kochanski, 2009) is

$$P(\mathrm{R}=i|\mathrm{S}) = \int_s^\infty \varphi(c_i;\mu_{C_i},\sigma_{C_i})\prod_{j\neq i}^M [1-\int_s^{c_i}\varphi(c_j;\mu_{C_j},\sigma_{C_j})dc_j]dc_i. \qquad (4)$$

*Computations for Maximum Likelihood Solutions*

   *Problem-specific module for BMCMC.* The problem-specific module for the Law of Categorical Judgment (Corrected) is lawcjcrk.cpp, and the overall program entry point is in bsr_analysis.py. Both are available in release z2009bsr_analysis-0.3.0.tar.gz. Each parameter can be free or constrained to some particular value.

   The module uses Romberg integration (available in the release gpklib-20090626.tar.gz and similar to Press et al., 2002, pp. 144-146) to compute a value of log*L*, the logarithm of the likelihood that the parameters would generate the actual rating frequencies. This value is returned to the BMCMC python modules. The problem-specific module also contains functions for input and calculates and writes all the required output. The latter includes the final parameters and the theoretical probabilities $P(R=i\,|\,S_h)$, their respective confidence limits, and values for measures of goodness of fit described below.

   Markov chain Monte Carlo is normally cast in Bayesian terms. Therefore, one must (formally) consider prior distributions. For this particular application, however, any reasonable prior will be much broader than the Bayesian posterior distribution, $\pi(\vec{p})$. Hence, the difference between the Bayesian and the likelihood approaches



becomes moot. One can readily interpret $L(\vec{p})$ as proportional to the (Bayesian) conditional probability of observing specified rating frequencies, given parameters $\vec{p}$. Samples from $\pi(\vec{p})$ will then reflect the distribution of parameters that are statistically consistent with the ratings.

*Parameter constraints and adjustments.* Before calculating $\log L$, the problem-specific module imposes four constraints that harmlessly remove degeneracies in the solutions. First, the mean of the signal density for $S_1$ is fixed at zero; all other signal means are made positive and are sorted into ascending order. Second, ascending order is imposed on criterion means (though they may be negative). Third, any negative signal or criterion standard deviation is made positive. Fourth, the average of all unconstrained signal and criterion standard deviations is divided into each unconstrained mean or standard deviation. The average standard deviation is thereby held at unity, keeping standard deviations in a convenient range. This rescaling is analogous to the conventional procedure of constraining a signal density standard deviation to unity. The constraints do not affect any results; they merely collapse together observationally equivalent solutions; and they make interpretation of the results easier. (These four are implemented via mcmc.problem_definition.fixer in mcmc.py; see the Appendix.)

A final adjustment effectively adds a Bayesian prior to $\log L$ to keep each standard deviation $\sigma$ away from zero. Gaussian densities take the form $\exp(-\ldots/\sigma^2)$. If $\sigma$ becomes small, say around 0.001, the integrands in Equation 1 change rapidly. In



turn, the Romberg routine forces the integration steps to be smaller than σ. A severe slowing of computation ensues. The program may even terminate prematurely if σ becomes so small that round-off errors prevent quadrature from achieving the desired accuracy. To avoid these difficulties, we add $\sum_u 0.1/\sigma_u$ to $\log L$, where $\sigma_u$ is any unconstrained standard deviation. When all $\sigma_u$ exceed 0.1, as generally happens here, this term hardly affects the solution.

Constraints and adjustments in optimisation mode must be employed carefully, however strong may be their pragmatic rationale. They can cause incorrect results if they have noticeable effects at or near a maximum likelihood solution. Badly formulated constraints also can artificially limit the size of the confidence intervals obtained in sampling mode. Nonetheless, adding smooth constraints is often much more efficient computationally than a hard-wall constraint (e.g. setting $\log L(\vec{p})$=-∞ outside the allowed region for $\log L(\vec{p})$).

*Goodness of fit.* In line with Schunn and Wallach (2005), we use four different ways of comparing the input rating frequencies against the  frequencies predicted from the recovered maximum likelihood parameters.  To begin with, we divide each input frequency by $Tr$/S$_h$ expresses it as a proportion  $p(R=i\,|\,S_h)$ of responses $i$ to stimulus $h$. Then, the theoretical conditional probabilities $P(R=i\,|\,S_h)$ are computed via Equations 1, 3, or 4 from the optimal parameters. The four procedures for evaluating goodness of fit require calculation of (a) root-mean-square deviation



(*RMSD*) between $p(R = i \mid S_h)$ and $P(R = i \mid S_h)$, (b) Pearson's $r^2$ for the regression of

$p(R = i \mid S_h)$ on $P(R = i \mid S_h)$, (c) the regression coefficients $b_0$ and $b_1$ and their 95 per

cent confidence limits and (d) Kullback-Leibler divergence, the relative entropy of

two probability densities $p$ and $q$ (Kullback & Leibler, 1951; Weisstein, n.d.),

specified here as $\sum_h \sum_i P(R = i \mid S_h) \log_2 [P(R = i \mid S_h) / p(R = i \mid S_h)]$.

*Tests of The BMCMC Procedure*

We used BMCMC to fit the FSDT, SDT, and CSDT models (Equations 1, 3,

and 4, respectively) to pseudo data matrices generated from known parameter values.

The parameters are the means and standard deviations of the densities in a FSDT,

SDT, or CSDT model. One parameter set was chosen independently for each type of

model. Signal and/or criterion standard deviations (as appropriate) were varied

irregularly, as were distances between signal and criterion means.

Pseudo data matrices were generated by trial-by-trial simulation of a rating

experiment with six stimuli and 10 responses (i.e. 9 criteria). On each trial, a sample

was drawn randomly from a signal density along with one sample from each of the 9

criterion densities. Using the decision rule, a response was selected. From each

parameter set, we produced three pseudo data matrices, with 200, 500, and 1000 trials

per stimulus (*Tr*/$S_h$), respectively.

Using the appropriate model equation, maximum likelihood parameters were

computed for each pseudo data matrix through the BMCMC algorithm. Ninety-five

per cent confidence limits for the recovered parameters and the theoretical



probabilities $P(R = i \mid S_h)$ were based on more than 4000 sets of parameters produced by BMCMC in sampling mode. Three separate fits were made to each pseudo data matrix, starting from different initial guesses at the parameters. The fit with the highest log likelihood (excluding the soft constraint $\sum_u 0.1 / \sigma_u$) was accepted. Altogether, we undertook 27 fits (3 models × 3 matrices × 3 starting points).

We assessed the consistency among the three fits to a given matrix. The absolute difference $\mid \log L_{\max} - \log L_{\min} \mid$ was found between the highest and the lowest of the three $\log L$ values. We divided this difference by the degrees of freedom in the pseudo data matrix. The result for fits of the FSDT, SDT, and CSDT models was always less than 0.05, indicating good consistency between the three fits to each pseudo data matrix. All optimisations on a given matrix apparently gave good approximations to the true maximum likelihood solution.

To compare the generating parameters against those recovered by an optimisation, any difference in scaling units had to be eliminated. We obtained a least-squares solution for $b$ in the equation $R = bG$, where $R$ and $G$ are the recovered and generating parameters, respectively. Then the recovered parameters $R$ and their confidence limits were plotted against the rescaled generating parameters $G_T = bG$. For a successful recovery, the points should fall on or near the major diagonal, and 95 per cent of the confidence limits should intersect that line.

*Recovery results*. Figure 1 shows the best set of recovered parameters plotted against the rescaled generating values for the FSDT-generated pseudo data matrices.



The upper, middle, and lower panels are plots for the matrices with 200 $Tr$/S$_h$ (FSDT200), 500 Tr/S$_h$ (FSDT500), and 1000 Tr/S$_h$ (FSDT1000), respectively. Ninety-five per cent confidence limits appear around each recovered parameter value.

_______________________

Figure 1 about here

_______________________

The left side of each panel displays four items: the rescaling equation for the generating parameters; log$L_G$, the log likelihood produced by evaluating the match of the pseudo data to the generating parameters; log$L_R$, the log likelihood for the recovered parameters; and the number of samples $N_e$ on which the confidence limits are based. The legend in Figure 1A applies to each panel.

The BMCMC algorithm always found a solution whose log$L_R$ exceeded the original log$L_G$. However surprising this may be at first, it should be expected. The pseudo data proportions generated by simulation differ randomly from the underlying probabilities yielded by the generating parameters. Consequently, parameters will typically occur that fit better than the generating parameters, and a successful optimisation should find them. Nonetheless, 95 per cent of the rescaled generating parameters should fall within the confidence limits for the recovered parameters.



The points in each panel of Figure 1 follow the major diagonal well. Naturally, confidence intervals shrink as $Tr/S_h$ increases. All confidence limits in the figure intersect the major diagonal. In short, the fits and error bars are perfectly satisfactory.

Figures 2 and 3 show the results for the SDT and CSDT fits, respectively. They are organized like Figure 1. Recovered $\log L_R$ was again higher than generating $\log L_G$. In each of Figures 2C, 3B, and 3C, confidence limits fail once to intersect the major diagonal. This is predictable: 5 per cent of the total 63 parameters should yield misses over the long run. The fits to the SDT and CSDT pseudo data matrices also are excellent.

---------------------------------

Figures 2 and 3 about here

---------------------------------

*Goodness of fit*. Table 1 presents information on goodness of fit for the results in Figures 1 through 3. The first three rows refer to the fits of the FSDT model to the FSDT-generated pseudo data matrices. For these fits, root mean square deviation (*RMSD*) was 0.014 at most. The linear regression of $p(R=i|S_h)$ on $P(R=i|S_h)$ yielded an $r^2$ of at least .95. The 95 per cent confidence limits for the regression coefficients $b_1$ and $b_0$ included the ideal values of unity and zero, respectively. Finally, the Kullback-Leibler divergence always fell well below 0.1.



All indicators of goodness of fit tended to improve as $Tr/S_h$ increased. The high $r^2$ for the 200 $Tr/S_h$ matrix actually left little room for such increases. The Kullback-Leibler values declined as $Tr/S_h$ increased, because the divergence acts as an error measure on differences between pseudo-data proportions and theoretical probabilities. These differences include statistical fluctuations that are larger for smaller experiments. In short, all calculations indicate good fits of the recovered theoretical probabilities to the pseudo data proportions.

The results in Table 1 for the CSDT and SDT intra-model fits repeat those for the FSDT fits. The $RMSD$ values varied between .005 and .012. Again, $r^2$ exceeded .95. The confidence limits for $b_1$ and $b_0$ always included unity and zero, respectively. Finally, Kullback-Leibler divergence values were small and decreased from about 0.06 for $Tr/S_h$=200 to about 0.01 for $Tr/S_h$=1000.

_______________________

Table 1 about here

_______________________

*Cross-Model Fits*

The widespread use of signal detection theory motivates an obvious question. For practical-size experiments, could SDT still fit data generated by the more complex FSDT model? To obtain an initial answer, we fitted the SDT model (Equation 3) to the pseudo data generated by the FSDT model (Equation 1) and



compared the results to those from the previous FSDT intra-model fits. We extended

this cross-model procedure by fitting the CSDT model to the FSDT-generated pseudo

data and fitting the SDT model to the CSDT-generated pseudo data. Again, three

different starting points were used for each fit. This gave nine SDT-to-FSDT fits (3

matrices × 3 starting points), nine CSDT-to-FSDT fits, and nine SDT-to-CSDT fits.

As before, the fits from the three different starting points proved consistent. For each

cross-model fit, however, $\log L_R$ was always somewhat lower than for the

corresponding intra-model fit.

    *Goodness of fit*. Table 2 contains the indicators of goodness of cross-model

fits. Here, *RMSD* always exceeded 0.02, twice the typical values for intra-model fits.

Furthermore, $r^2$ was .949 at best, even falling to .878. The confidence limits for the

regression coefficients $b_0$ and $b_1$ always exceeded those from the FSDT and CSDT

intra-model fits. Those limits, however, still included the theoretically desirable

values of zero and unity for $b_0$ and $b_1$, respectively. Most noticeably, Kullback-Leibler

divergence always exceeded 0.1, with a median of 0.23.  This is above the typical

value for the intra-model fits with $Tr/S_h$=200 and is about 20 times larger than the

values for intra-model fits with $Tr/S_h$=1000.

       The findings in Table 2 suggest that Kullback-Leibler divergence is a sensitive

indicator of goodness of fit of model probabilities to observed proportions.

Confidence limits for the regression coefficients and *RMSD* can also be useful.

Splitting the differences between Tables 1 and 2 leads to the following rough guides



to acceptable fits for experiments of the size studied here: *RMSD* should be less than

.025; confidence limits for $b_0$ and $b_1$ should remain below ±0.0075 and ±0.05,

respectively; and Kullback-Leibler divergence should be less than 0.175. Under these

guidelines, the SDT fit to the FSDT-500 matrix would be considered marginal. No

other cross-model fit would be acceptable.

________________________

Table 2 about here

________________________

      The cross-model results indicate that a subject who obeys the FSDT rating

model may produce data beyond the reach of the classical SDT model. Indeed, that

model cannot entirely explain the findings of one study on hearing. Pastore and

Macmillan (2002) undertook an SDT reanalysis of eight sets of rating data originally

collected by Schouten and van Hessen (1998) for an intensity and for a speech

continuum. Except for one data set for intensity, Pastore and Macmillan found that

the SDT model fitted only some ROC curves. Luce's (1963) low threshold model

seemed to account for the numerous exceptions. This left behind a situation where

two different models had to be invoked to explain single data sets. The FSDT model

may provide an alternative to this unsatisfactory outcome.



CONCLUSIONS

In situations requiring numerical optimisation, the Bootstrap Markov Chain Monte Carlo algorithm can efficiently provide maximum likelihood solutions and confidence limits for the parameters. Tests of the algorithm gave good fits of each variant of the Law of Categorical Judgment (Corrected) to pseudo data matrices generated by that same variant. Furthermore, measures of goodness of fit showed that these intra-model fits were better than fits of inappropriate models to the pseudo data. The way is now open to complete analyses of experimental rating data with the Gaussian based Law of Categorical Judgment (Corrected).



Appendix

Hard-Wall Constraints and Symmetries

Multiple hard-wall constraints can slow the BMCMC algorithm. However, hard-wall constraints can often be eliminated by transforming the problem into a symmetrised version. For instance, consider a truncated normal density $\pi(y) = \{ n(y; \mu, \sigma)$ if $y > 0$; else $0 \}$. With a hard wall at $y = 0$, samples can be obtained from a symmetrised variant of $\pi$, $\zeta(q) \propto n(|q|; \mu, \sigma)$, and then mapping $y = |q|$. (Note that $\zeta(q)$ is bimodal if $\mu > 0$.) This can be implemented by setting up the "fixer" method in a class derived from problem_definition so that it maps negative $y$ into positive $y$ at each step:

```
class problem_definition(object):
        def fixer(self, y):

                return numpy.absolute(y)
```

This kind of substitution can dramatically improve speed, especially if there are multiple constraints. However, two conditions must be met. First, the mapping needs to be (at least locally) a reflection around a plane. Second, this technique is known to work only if the mapping $y = M(q)$ is such that $\pi(M(q)) = \zeta(q)$ for all $q$, and $M$ is piecewise isometric linear.

In practice, this works out conveniently if the constraints are orthogonal or parallel hyperplanes. For instance, if we have a vector $\vec{p} = (x, y)$, the constraints $0 < x < 1$ and $y > 0$ are easy to implement.

Author Note

Greg Kochanski and Burton S. Rosner, Phonetics Laboratory, University of Oxford.



Footnotes

[1]This has yielded good behavior in certain test cases, but no proofs are available. To do this, use mcmc.position_nonrepeatable instead of the default mcmc.position_repeatable.

[2]In the code, this can be done by implementing mcmc.problem_definition.fixer.

[3]By means of mcmc.problem_definition.logP.

[4]See the BootStepper class in mcmc.py.

[5]See the StepV method in the BootStepper class in mcmc.py.

[6]See mcmc.adjuster.inctry_guts and mcmc.BootStepper.step_boot.

[7]See _inctry_guts in the adjuster class in mcmc.py.

[8]The algorithm is reasonably robust to a mismatch between V and $\pi$. In practice, it frequently tolerates a two orders of magnitude mismatch in standard deviation, especially in optimization mode. Starting in sampling mode can be slow if the initial step acceptance probability is low.

[9]The high level interface to optimization mode is mcmc_helper.stepper.run_to_bottom.

[10]The annealing schedule can be redesigned by changing parameters or re-implementing the step_acceptor object in mcmc_helper.py.

[11]See mcmc.BootStepper.reset.

[12]See code that uses mcmc.Archive.sorted.



[13]Optimization mode is the run_to_bottom method in the stepper class in mcmc_helper.py. See the _not_at_bottom method of the stepper class in mcmc_helper.py.

[14]Direction changes are accumulated in the dotchanged and xchanged variables in stepper.run_to_bottom in mcmc_helper.py.

[15]Computed in Bootstepper.ergodic in mcmc.py.

[16]This information is accumulated in the "es" variable in stepper.run_to_bottom in mcmc_helper.py. N.B.: if there are $k \geq 2$ comparably large eigenvalues of the covariance matrix, this termination condition may fire early; it assumes that there is a single long axis to $\pi$.

[17]The high-level interface to sampling mode is mcmc_helper.stepper.run_to_ergodic.

[18]See the "nc" variable in stepper.run_to_ergodic in mcmc.py.



Table 1

*Goodness of Intra-Model Fits*

_______________________________________________________________

| Model[a] | $Tr/S_h$[b] | Indicator | | | | |
|---|---|---|---|---|---|---|
| | | $RMSD$[c] | $r^2$ | $b_0 \pm CL95$[d] | $b_1 \pm CL95$ | $K\text{-}B$[e] |
| FSDT | 200 | .014 | .953 | 0.00±0.005 | 1.014±0.031 | 0.059 |
| | 500 | .009 | .960 | 0.000±0.003 | 1.001±0.021 | 0.028 |
| | 1000 | .005 | .965 | 0.000±0.002 | 1.004±0.012 | 0.009 |
| SDT | 200 | .012 | .960 | 0.001±0.004 | 0.991±0.022 | 0.062 |
| | 500 | .009 | .963 | 0.001±0.003 | 0.999±0.016 | 0.024 |
| | 1000 | .007 | .965 | 0.000±0.002 | 1.001±0.013 | 0.014 |
| CSDT | 200 | .010 | .964 | 0.000±0.003 | 1.001±0.015 | 0.046 |
| | 500 | .006 | .966 | 0.000±0.002 | 1.005±0.010 | 0.014 |
| | 1000 | .005 | .966 | 0.000±0.001 | 1.001±0.007 | 0.010 |

[a]See text for model type designations. [b]Trials per stimulus. [c]Root mean square difference between pseudo proportions and recovered probabilities of rating $S_h$ as *i*. [d]Double-sided 95 per cent confidence limits. [e]Kullback-Leibler divergence.



Table 2

*Goodness of Cross-Model Fits*

| Model[a] | | $Tr$/S$_h$[b] | Indicator | | | | |
|---|---|---|---|---|---|---|---|
| Generating | Fitted | | $RMSD$[c] | $r^2$ | $b_0 \pm CL95$[d] | $b_1 \pm CL95$ | $K\text{-}B$[e] |
| FSDT | SDT | 200 | .030 | .904 | -0.005±0.010 | 1.045±0.072 | 0.230 |
| | | 500 | .022 | .932 | -0.001±0.008 | 1.011±0.051 | 0.143 |
| | | 1000 | .024 | .925 | -0.004±0.008 | 1.038±0.058 | 0.195 |
| | CSDT | 200 | .038 | .868 | -0.006±0.014 | 1.064±0.093 | 0.405 |
| | | 500 | .035 | .875 | -0.005±0.013 | 1.051±0.088 | 0.351 |
| | | 1000 | .033 | .890 | -0.006±0.012 | 1.059±0.080 | 0.316 |
| CSDT | SDT | 200 | .029 | .937 | 0.000±0.009 | 1.001±0.046 | 0.239 |
| | | 500 | .023 | .948 | -0.001±0.007 | 1.006±0.037 | 0.207 |
| | | 1000 | .022 | .949 | 0.000±0.007 | 1.002±0.036 | 0.192 |

[a]See text for model type designations. [b]Trials per stimulus. [c]Root mean square difference between pseudo proportions and recovered probabilities of rating S$_h$ as *i*. [d]Double-sided 95 per cent confidence limits. [e]Kullback-Leibler divergence.



Figure Captions

*Figure 1*. Recovered and generating parameters for 6 by 10 matricies of rating pseudo

data generated by the FSDT model. Confidence limits are 95 per cent. A. Matrix from

trial-by-trial simulation, 200 trials per stimulus (*Tr*/S$_h$). B. Simulated matrix, 500

*Tr*/S$_h$. C. Simulated matrix, 1000 *Tr*/S$_h$. See text for further explanation.

*Figure 2*. Same as Figure 1, but for the SDT model.

*Figure 3*. Same as Figure 1, but for the CSDT model.



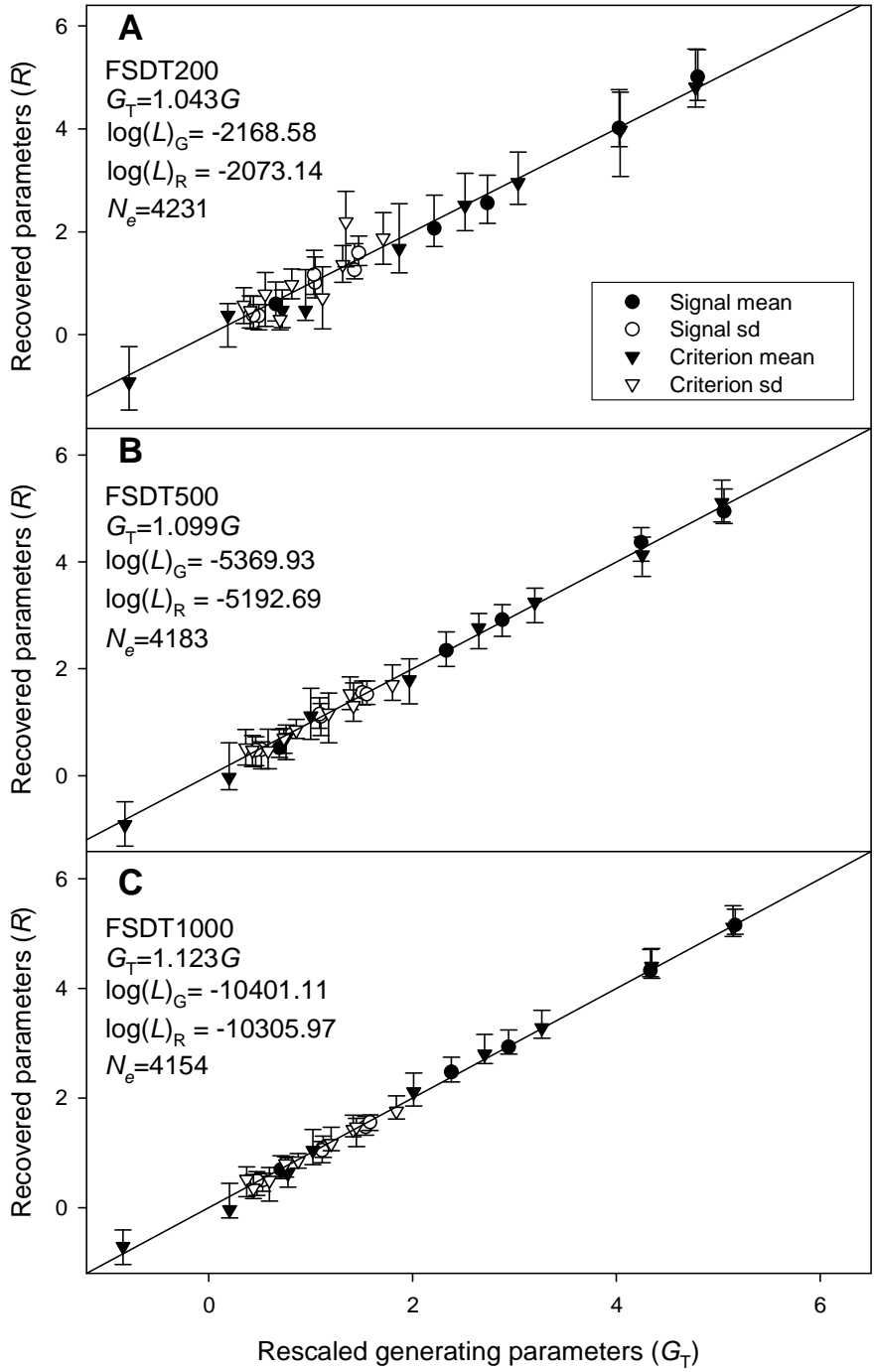



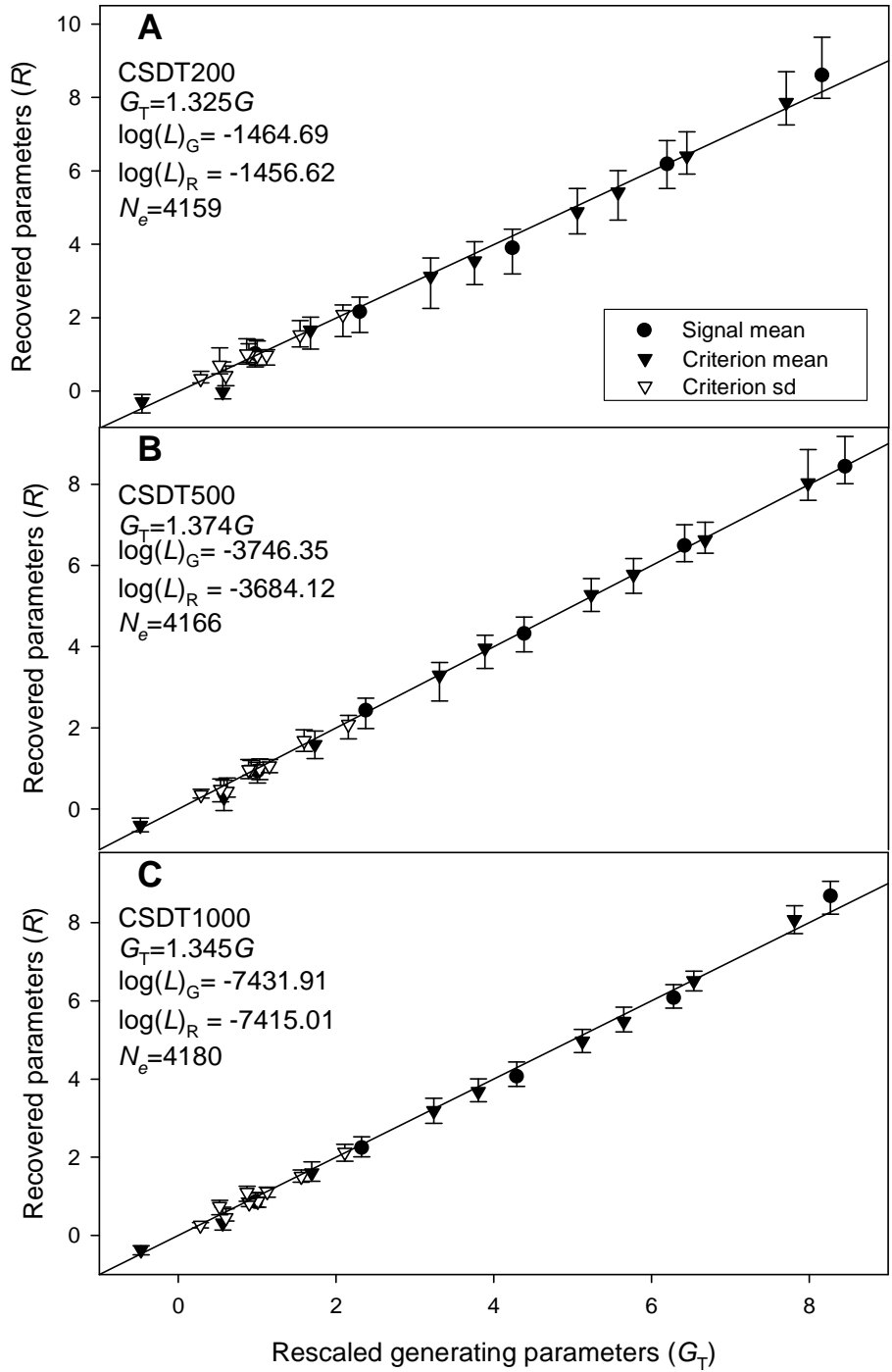



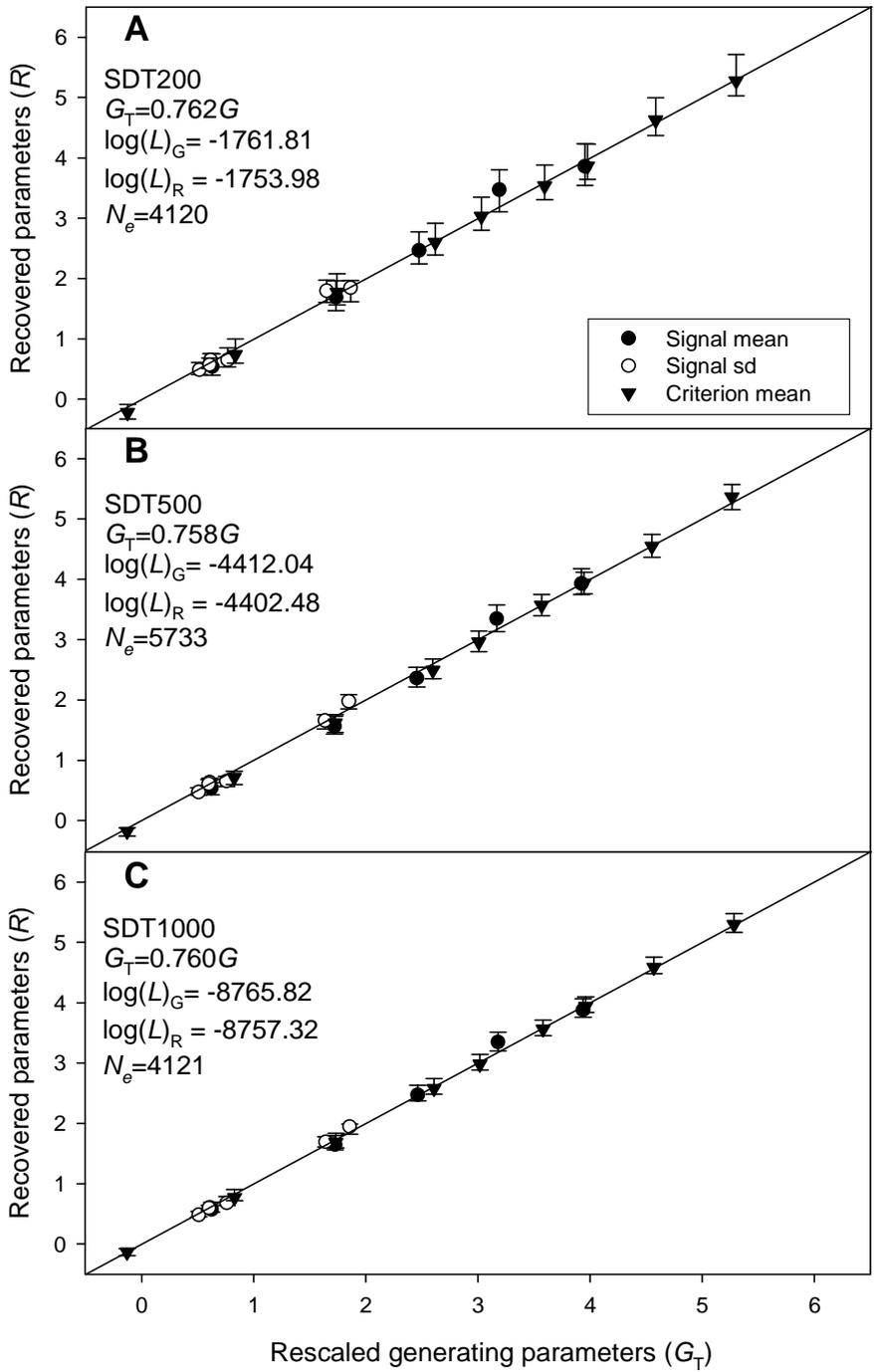